\documentclass[12pt,reqno]{amsart}
\setbox0=\hbox{$+$}
\newdimen\plusheight
\plusheight=\ht0
\def\+{\;\lower\plusheight\hbox{$+$}\;}

\setbox0=\hbox{$-$}
\newdimen\minusheight
\minusheight=\ht0
\def\-{\;\lower\minusheight\hbox{$-$}\;}

\setbox0=\hbox{$\cdots$}
\newdimen\cdotsheight
\cdotsheight=\plusheight%\ht0
\def\cds{\lower\cdotsheight\hbox{$\cdots$}}

\renewcommand{\(}{\left\(}
\renewcommand{\)}{\right\)}

\def\NI{\noindent}

\renewcommand{\pmod}[1]{\,(\textup{mod}\,#1)}
\numberwithin{equation}{section}
\theoremstyle{plain}
\newtheorem{theorem}{Theorem}[section]
\newtheorem{lemma}[theorem]{Lemma}

\begin{document}
\begin{center}{\bf \large  Some New Congruences For Overpartition Function With $\ell$-Regular Non-Overlined Parts}\end{center}\vskip
 5mm
\begin{center}
\centerline{\bf Nipen Saikia$^{\ast}$ and Adam Paksok}\vskip2mm

{\it Department of Mathematics, Rajiv Gandhi University,\\ Rono Hills, Doimukh, Arunachal Pradesh, India-791112.\\ E. Mail(s): nipen.saikia@rgu.ac.in; adam.paksok@rgu.ac.in}\\
 	$^\ast$\textit{Corresponding author}\end{center}\vskip3mm
% \footnotetext[1]{Corresponding author.}

 \noindent{\bf Abstract:} Alanzi et al. (2022) investigated overpartition of a positive integer $n$ with $\ell$-regular non-overlined parts denoted by $\overline R_\ell^\ast (n)$, and proved some results for the case $\ell=3$. As extension to the results of Alanzi et al., Sellers (2024) proved some new congruences for $\overline R_3^\ast (n)$.  In this paper, we prove some new infinite families and particular congruences for $\overline R_\ell^\ast (n)$ for $\ell=4, 5k, 6$, and  8, where $k$ is any positive integer. We also offer some congruences connecting $\overline R_\ell^\ast (n)$ with some other partition functions.

 \vskip
 3mm

 \vskip
 3mm

 \NI{\bf Mathematics Subject Classifications:} 11P83;  05A17. \vskip
 3mm

 \section{Introduction}
  A partition of a positive integer $n$ is a finite non-increasing sequence of positive integers $s_1, s_2, \dots, s_k $ satisfying $s_1+s_2+\cdots +s_k=n$, where the integers $s_1, s_2, \dots, s_k $  are called as parts or summands of the partition. For example, the partitions of $n=3$ are given by $3, 2+1, 1+1+1$. The number of partition of a non-negative integer $n$ is usually  denoted by $p(n)$ and its  generating function is given by 
  \begin{equation}\label{pngeni}\sum_{n=0}^{\infty}p(n)q^n=\dfrac{1}{(q;q)_{\infty}},\qquad p(0)=1,\end{equation} where for any complex number $a$ and $|q|<1$, 
  \begin{equation}\label{9}
  (a;q)_{\infty}=\prod_{n=0}(1-aq^n).\end{equation} For brevity, we will use the notation  
  $f_h:=(q^h;q^h)_{\infty}$ for any positive integer $h$. Ramanujan \cite{sr2, SR}(also see \cite{bc1}) offered following beautiful congruences for the partition function $p(n)$:
  \begin{equation}
  \label{pc}
  p(5n+4) \equiv0\pmod5,\quad
  p(7n+5)\equiv0 \pmod7,\quad p(11n+6)\equiv0\pmod {11}.\end{equation}  
  For any positive integer ${\ell}$, an ${\ell}$-regular partition of a positive integer $n$ is a partition of $n$ such that none of its part is divisible by ${\ell}$.  If $b_{\ell}(n)$ denotes the number of ${\ell}$-regular partitions of $n$, then the generating function of $b_{\ell}(n)$ is given by  
   \begin{equation}\label{lregularintro}
    	\sum_{n=0}^{\infty} b_{\ell}(n)q^n=\frac{(q^{\ell};q^{\ell})_\infty}{(q;q)_\infty}=\frac{f_\ell}{f_1}.\end{equation}
   The arithmetic properties of ${\ell}$-regular partition have been studied by different authors, for example see \cite{zlndb, klndb, carlson, cui, penniston, hirs, lin2, penniston1, webb, xia1}.
 
 An overpartition of a positive integer $n$ is a partition of $n$ wherein first occurrence of any part may be overlined.  The number of overpartitions of a positive integer $n$ is usually denoted by $\overline p(n)$ which is first studied by Corteel and Lovejoy \cite{cor}. The generating function of $\overline p(n)$ \cite{cor} is given by
 \begin{equation}\label{opn}\sum_{n\ge0}^{\infty}\overline p(n)q^n=\dfrac{(q^2;q^2)_\infty}{(q;q)_{\infty}^2}=\dfrac{f_2}{f_1^2}.\end{equation}  
 For instance, $\overline p(3)=8$ with overpartitions of 3 given by  $3, \overline{3}, 2+1, \overline 2+1, 2+\overline 1, \overline 2+\overline 1, 1+1+1, \overline1+1+1$.
 One may see \cite{hs1, hs2, kim1, kmah, ij, xia1over} and references therein for congruence and other arithmetic properties of overpartition function $\overline{p}(n)$. 
 
 As an extension to the overpartition function $\overline p(n)$, overpartition with $\ell$-regular overlined parts  denoted by $\overline p_\ell(n)$ was defined. For example, $\overline p_2(3)=6$ with the relevant partitions given by
  $ 3, \; \overline{3}, \;2+1, \; 2+\overline{1},\; 1+1+1,\; \overline{1}+1+1$.
   The generating function for $\overline p_\ell(n)$ \cite{nak} is given by
   \begin{equation}\label{cq1}\sum_{n=0}^{\infty} \overline p_\ell(n)q^n=\frac{(q^2;q^2)_\infty (q^{\ell};q^{\ell})_\infty}{(q;q)_\infty^2(q^{2\ell};q^{2\ell})_\infty}=\dfrac{f_2f_\ell}{f_1^2f_{2\ell}}.
     \end{equation}
   Ballantine and Merca \cite{nak} investigated the connections between $\overline p_\ell(n)$ and  other restricted overpartition functions, and also established recurrence relations for $\overline p_3(n)$.  Shivanna and Chandrappa \cite{sc} proved congruences modulo $5, 7, 40 $ and powers of 2 for  $\overline p_\ell(n)$ with $\ell\in\{4,8\}$. 
 Buragohain and Saikia \cite{ps} proved some congruences for  $\overline p_\ell(n)$ with $\ell \in {3,5,9,12,16,18,24}$ using $q$-series manipulations.
 
 Recently, Alanzi et al. \cite{alanzi} considered overpartition of positive integer $n$ with $\ell$-regular non-overlined parts denoted by $\overline R_\ell^\ast (n)$. For instance, $\overline R_2^\ast(3)=6$ with the relevant partitions given by  $ 3, \; \overline{3}, \;\overline 2+1,\; \overline 2+\overline{1}\; 1+1+1\; \overline{1}+1+1$.
Alanzi et al. \cite{alanzi} proved some general parity results and  specific congruences for the case $\ell=3$  of $\overline R_\ell^\ast (n)$. The generating function of $\overline R_\ell^\ast (n)$ was given by Sellers \cite{sq} as
 \begin{equation}\label{rast}
 \sum_{n=0}^{\infty}  \overline R_\ell^\ast (n)q^n=\frac{ (q^{\ell};q^{\ell})_\infty(q^2;q^2)_\infty}{(q;q)_\infty^2}=\dfrac{f_2f_\ell}{f_1^2}.
 \end{equation} As an extension to the work of Alanzi et al. \cite{alanzi},  Sellers \cite{sq} found some infinite families and some specific congruences modulo $3$ and $4$ for  $\overline R_3^\ast (n)$. In closing remarks of \cite{sq}, further congruences for the partition function $\overline R_\ell^\ast (n)$ are desired. 
 
 In this paper, we will prove some infinite families and particular congruences for $\overline R_\ell^\ast (n)$ for $\ell=4, 5k, 6, 8 $, where $k$ is any positive integer. We also offer some congrunce relations connecting $\overline R_\ell^\ast (n)$ with some other partition functions.

  \section{Preliminaries}
   Ramanujan's general theta-function $F\left(x,y\right)$ is defined by
   \begin{equation}\label{eq2a}
    F(x, y)=\sum_{n=-\infty}^\infty x^{n(n+1)/2}y^{n(n-1)/2},~  |xy|<1.
   \end{equation} 
   Three important special  cases of $F\left(x,y\right)$ \cite[p. 36, Entry 22 (i), (ii))]{bcb3} are  the theta-functions $\phi(q)$, $\psi(q)$  and  $f(-q)$ given by
   \begin{equation}\label{e0}
  \phi(q):=F(q, q)=\sum_{\nu=0}^{\infty}q^{\nu^2}=(-q;q^2)^{2}_{\infty}(q^2;q^2)_{\infty}=\dfrac{{f
  _2^5}}{{f_1^2}{f_4^2}} ,
   \end{equation}
  \begin{equation}\label{tp1}
 \hspace{-1cm} \psi(q):=F(q,q^3)=\sum_{\nu\geq 0}q^{\nu(\nu+1)/{2}}=\dfrac{(q^2;q^2)_{\infty}}{(q;q^2)_{\infty}}=\dfrac{{f
  _2^2}}{f_1},
  \end{equation}
   and
  \begin{equation}\label{t3}
  \hspace{-.7cm}	f(-q):=F(-q,-q^2)=\sum_{\nu=-\infty}^{\infty}(-1)^\nu q^{\nu(3\nu+1)/{2}}=f_1.\end{equation}
   Employing elementary $q$-operations, it is easily seen that \begin{equation}\label{phim}
  \phi(-q)=\dfrac{f_{1}^2}{f_2}.
  \end{equation}
  
   \begin{lemma}\label{iii}\cite[Theorem 2.1]{cui} For any  odd prime, we have
  \begin{equation}\label{w4}
  \psi(q)= \sum_{j=0}^{(p-3)/2}{\it q}^{(j^2+j)/2} F\left( {\it q}^{(p^2+(2j+1)p)/2},{\it q}^{(p^2-(2j+1)p)/2}\right)+{\it q}^{(p^2-1)/8}\psi(q^{p^2}).\end{equation} \\Furthermore, $(j^2+j)/2 \not\equiv(p^2-1)/8 \pmod p~ for ~0 \leq j\leq (p-3)/2.$
  \end{lemma}
   
  \begin{lemma}\cite[Theorem 2.2]{cui} For any prime $p\geq 5$, we have
  $$ 
  f_1=\sum_{\substack{ k={-(p-1)/2} \\ k \ne {(\pm p-1)/6}}}^{(p-1)/2}(-1)^{k}{q}^{(3k^2+k)/2} F\left( {-q}^{(3p^2+(6k+1)p)/2},{-q}^{(3p^2-(6k+1)p)/2}\right)$$ \begin{equation}\label{u10}+(-1)^{{(\pm p-1)/6}}{q}^{(p^2-1)/{24}}f_{p^2},\end{equation}where
  \begin{equation*}
  \dfrac{\pm p-1}{6}
  = \left\{
          \begin{array}{ll} 
           \dfrac{(p-1)}{6},   
             & if~ p \equiv 1\pmod 6 \\
             \dfrac{(-p-1)}{6} ,
             & if~ p \equiv -1\pmod 6. 
          \end{array}
      \right.
  \end{equation*} \\Furthermore, if $\dfrac{-(p-1)}{2}\leq k \leq\dfrac{(p-1)}{2}~and ~k \neq \dfrac{(\pm p-1)}{6}$~ then\\
  $$\dfrac{3k^2+k}{2}\not\equiv \dfrac{p^2-1}{24} \pmod p.$$
  \end{lemma}
 
  \begin{lemma}\cite{ovyt} We have
  \begin{equation}\label{t1a}
  f_1^2=\dfrac{f_2{f_8^5}}{{f_4^2}{f_{16}^2}}-2q\dfrac{{f_2}f_{16}^2}{f_8},\end{equation}
  \begin{equation}\label{b4}
  \dfrac{1}{{f_1^2}}=\dfrac{{f_8^5}}{{f_2^5}{f_{16}^2}}+2q\dfrac{{f_4^2}{f_{16}^{2}}}{{f_2^5}{f_8}},\end{equation}
  \begin{equation}\label{b3}
  \dfrac{1}{{f_1^4}}=\dfrac{{f_4^{14}}}{{f_2^{14}}{f_8^4}}+4q\dfrac{{f_4^2}{f_8^4}}{{f_2^{10}}}
  \end{equation}
  and
  \begin{equation}\label{b4a}
  {f_1^4}=\dfrac{{f_4^{10}}}{{f_2^2}{f_8^4}}- 4q\dfrac{{f_2^2}{f_8^4}}{{f_4^2}}.
  \end{equation}
  \end{lemma}
  
  \begin{lemma}\cite{hs1} We have
  \begin{equation}\label{p1} 
  \dfrac{1}{\phi(-q)}=\dfrac{1}{\phi(-q^4)^4}\Big(\phi(q^4)^3+2q\phi(q^4)^2\psi(q^8)+4q^2\phi(q^4
  )\psi(q^8)^2+8q^3\psi(q^8)^3\Big).
  \end{equation}
  \end{lemma} 
  
    \begin{lemma}{\cite{hs1}}\label{le2} We have
    	\begin{align*}
    		\dfrac{1}{\phi(q)}&=\dfrac{\phi(q^{25})}{\phi(q^5)^6}\bigg[\phi(q^{25})^4-2q\phi(q^{25})^3X(q^5)+4q^2\phi(q^{25})^2X(q^5)^2-8q^3\phi(q^{25})X(q^5)^3\nonumber\\
    		&+q^4\left(16X(q^5)^4
    		-2\phi(q^{25})^3Y(q^5)\right)-12q^5\phi(q^{25})^2X(q^5)Y(q^5)+16q^6\phi(q^{25})X(q^5)^2Y(q^5)\nonumber\\
    		&-16q^7X(q^5)^3Y(q^5)+4q^8\phi(q^{25})^2Y(q^5)^2+16q^9\phi(q^{25})X(q^5)Y(q^5)^2\nonumber\\&+16q^{10}X(q^5)^2Y(q^5)^2
    		-8q^{12}\phi(q^{25})Y(q^5)^3-16q^{13}X(q^5)Y(q^5)^3+16q^{16}Y(q^5)^4\bigg],
    	\end{align*}
    	where 
    	\begin{equation*}\label{hF}
    		X(q)=\sum_{r=-\infty}^{\infty}q^{5r^2+2r} \quad \text{and}\quad Y(q)=\sum_{r=-\infty}^{\infty}q^{5r^2+4r}.
    	\end{equation*}
    \end{lemma}
  
  \begin{lemma}\cite[p. 49]{bcb3}\label{psii}
  We have
 $$\psi(q)=F(q^3,q^6)+q\psi(q^9).$$
  \end{lemma}
  
  We will alos employ following congruences which can be easily proved by using the binomial theorem: 
   For any prime $p$, 
  \begin{equation}\label{yp}
  {f_{p}}\equiv {f_1}^p \pmod p,\end{equation}
  and 
  \begin{equation}\label{yp2}
  {f_1^{p^2}}\equiv {f_p}^p \pmod {p^2}.\end{equation}

\section{Congrurences for $\overline R_\ell^\ast (n)$}

\begin{theorem}We have
$$\hspace{-8.8cm}(i)\quad\overline R_{4}^\ast (4n+\xi)\equiv 0\pmod 4; \quad \xi=2, 3,$$
$(ii)$\quad Let $p\ge 13$  be a prime with $\left(\dfrac{-6}{p}\right)=-1$ and $1 \le r \le p-1$. Then for any integers $\alpha\geq0$ and $n\geq0$, we have
\begin{equation}\label{q1}\sum_{n=0}^{\infty}\overline R_{4}^\ast\left(4\cdot p^{2\alpha}n+\frac{7\cdot p^{2\alpha}-1}{6}\right)\equiv 2f_1\psi(q^2)\pmod 4\end{equation}and
\begin{equation}\label{q2}\quad \overline R_4^\ast\left(4\cdot p^{2\alpha+1}(pn+r)+\frac{7\cdot p^{2\alpha+2}-1}{6}\right) \equiv 0\pmod4.\end{equation}
\end{theorem}
\begin{proof}Setting $\ell=4$ in \eqref{rast} and simplifying, we obtain
\begin{equation}\label{rast41}
	\sum_{n=0}^{\infty}  \overline R_{4}^\ast (n)q^n=\frac{ f_4f_2f_1^2}{f_1^4}.
\end{equation} 
Employing \eqref{yp2} with $p=2$ in \eqref{rast41} and then simplifying using \eqref{phim}, we obtain 
\begin{equation}\label{rast41a}
	\sum_{n=0}^{\infty}  \overline R_{4}^\ast (n)q^n\equiv f_4 \phi(-q)\pmod 4.
\end{equation}
From \cite[p.49]{bcb3}, we note that
\begin{equation}\label{4}
	\phi(q)=\phi(q^{n^2})+\sum_{r=1}^{n-1}q^{r^2}F(q^{n(n-2r)}, q^{n(n+2r)}).
\end{equation}
Setting $n=2$ in \eqref{4} and simplifying, we obtain
\begin{equation}\label{4a}
	\phi(q)=\phi(q^4)+q\psi(q^8).
\end{equation}
Replacing $q$ by $-q$ in \eqref{4a} and then employing in \eqref{rast41a}, we obtain
\begin{equation}\label{rast41b}
	\sum_{n=0}^{\infty}  \overline R_{4}^\ast (n)q^n\equiv f_4\left(\phi(q^4)-2q\psi(q^8)\right)\pmod{4}.
\end{equation}
Extracting the terms involving $q^{4n+\xi}$ for $\xi=2,3$ from \eqref{rast41b}, we arrive at (i). 

To prove (ii), we employ the principle of mathematical induction. Extracting the terms involving $q^{4n+1}$ from \eqref{rast41b}, dividing by $q$ and replacing $q^4$ by $q$, we obtain
\begin{equation}\label{rast41c}
		\sum_{n=0}^{\infty}  \overline R_{4}^\ast (4n+1)q^n\equiv2f_1\psi(q^2)\pmod{4},
\end{equation}
	which is $\alpha=0$ case of (ii). Assume that (ii) is true for some integer $\alpha\geq0$. Using \eqref{w4} and \eqref{u10} in (ii), we obtain
$$\hspace{-10cm}\sum_{n=0}^{\infty}\overline R_{4}^\ast\left(4\cdot p^{2\alpha}n+\frac{7\cdot p^{2\alpha}-1}{6}\right) $$
	$$\equiv2\bigg[\sum_{\substack{ k={-(p-1)/2} \\ k \ne {(\pm p-1)/6}}}^{(p-1)/2}(-1)^{k}{q}^{\frac{3k^2+k}{2}} F\left( {-q}^{(3p^2+(6k+1)p)/2},{-q}^{(3p^2-(6k+1)p)/2}\right)+(-1)^{{(\pm p-1)/6}}{q}^{(p^2-1)/24}f_{p^2}\bigg]$$
		\begin{equation}\label{rast41d}
		\hspace{1cm}\times\bigg[\sum_{j=0}^{(p-3)/2}{\it q}^{2(j^2+j)/2} F\left( {\it q}^{2(p^2+(2j+1)p)/2},{\it q}^{2(p^2-(2j+1)p)/2}\right)+{\it q}^{2(p^2-1)/8}\psi(q^{2p^2})\bigg]\pmod{4}.
	\end{equation}
	Consider the congruence
	$$\left(\frac{3k^2+k}{2}\right)+2\left(\frac{j^2+j}{2}\right)\equiv \frac{7(p^2-1)}{24}\pmod{p} ,$$
	which is similar to
	$$(6k+1)^2+6(2j+1)^2\equiv0\pmod{p}.$$
	Since $\left(\dfrac{-6}{p}\right)=-1$, the above congruence has only solution $k=(\pm p-1)/6$ and $j=(\pm p-1)/2$. So, extracting the terms involving $q^{pn+7(p^2-1)/24}$ from both sides of \eqref{rast41d}, dividing  by $q^{7(p^2-1)/24}$ and then replacing $q^p$ by $q$, we obtain
		\begin{equation}\label{rast41e}
		\sum_{n=0}^{\infty}\overline R_{4}^\ast\left(4p^{2\alpha+1}n+\frac{7p^{2\alpha+2}-1}{6}\right)q^n\equiv2f_p\psi(q^{2p})\pmod{4}.
	\end{equation}
	Extracting the terms involving $q^{pn}$ from \eqref{rast41e} and replacing $q^p$ by $q$, we obtain
	\begin{equation}\label{rast41f}
		\sum_{n=0}^{\infty}\overline R_{4}^\ast\left(4p^{2(\alpha+1)}n+\frac{7p^{2\alpha+2}-1}{6}\right)q^n\equiv2f_1\psi(q^2)\pmod{4},
	\end{equation}
	which is $\alpha+1$ case of \eqref{q1}. 	
	Extracting the coefficients of the terms involving $q^{pn+r}$, for $1\le r\le p-1$, from both sides \eqref{rast41e}, we arrive \eqref{q2}.
\end{proof}

\begin{theorem}For any positive integer $k$, we have
$$\hspace{-8.5cm}(i)\quad\overline R_{5k}^\ast (5n+\xi)\equiv 0\pmod 4; \quad \xi=2, 3,$$
$$\hspace{-10.5cm}(ii)\quad\overline R_{5k}^\ast (5n+1)\equiv 0\pmod 2.$$
\end{theorem}
\begin{proof}Setting $\ell=5k$ in \eqref{rast} and employing \eqref{phim}, we obtain
	\begin{equation}\label{rast5k}
	\sum_{n=0}^{\infty}  \overline R_{5k}^\ast (n)q^n=~\frac{f_{5k}f_2}{f_1^2}=\frac{ f_{5k}}{\phi(-q)}.
	\end{equation} Employing  Lemma \ref{le2} in \eqref{rast5k} and simplifying using \eqref{yp2} with $p=2$, we obtain
$$\sum_{n=0}^{\infty}  \overline R_{5k}^\ast (n)q^n\equiv  \dfrac{f_{5k}~ \phi(-q^{25})}{\phi(-q^5)^6}\Big[\phi(-q^{25})^4+2q\phi(-q^{25})^3X(q^5)$$
\begin{equation}\label{po1}-2q^4
	\phi(-q^{25})^3Y(-q^5)\Big]\pmod 4. \end{equation}
	Extracting the terms involving $q^{5n+\xi}$ for $\xi=2, 3$ from  \eqref{po1}, we arrive at (i).  Similarly, (ii) follows from \eqref{rast5k} by employing \eqref{yp} with $p=2$ and then extracting the terms involving $q^{5n+1}$.
\end{proof}

\begin{theorem} For any integer $\alpha\ge 0$, we have
$$\hspace{-7.9cm}(i)\quad \overline R_6^\ast\left(9^{\alpha}n +\frac{9^\alpha-1}{4}\right) \equiv \overline R_6^\ast(n)\pmod3,$$
$$\hspace{-7.6cm}(ii)\quad \overline R_6^\ast\left(9^{\alpha+1}n +\frac{21\cdot 9^\alpha-1}{4}\right) \equiv 0\pmod3 ,$$
$$\hspace{-7.6cm}(iii)\quad \overline R_6^\ast\left(9^{\alpha+1}n +\frac{33\cdot 9^\alpha-1}{4}\right) \equiv 0\pmod3.$$
\end{theorem}
\begin{proof} Setting $\ell=6$ in \eqref{rast}, we obtain
	\begin{equation}\label{ras6}
	\sum_{n=0}^{\infty}\overline R_{6}^\ast (n)q^n=\frac{f_6f_2}{f_1^2}.\end{equation}
Simplifying \eqref{ras6} by employing \eqref{yp} with $p=3$ and, we obtain
\begin{equation}\label{ras61}
\sum_{n=0}^{\infty}\overline R_{6}^\ast (n)q^n\equiv\frac{f^4_2}{f_1^2}\pmod 3.\end{equation}
Employing  \eqref{tp1} in \eqref{ras61}, we obtain
\begin{equation}\label{ras61a}
	\sum_{n=0}^{\infty}\overline R_{6}^\ast (n)q^n\equiv\psi^2(q)\pmod 3.
\end{equation}
Employing Lemma \ref{psii} in \eqref{ras61a}, we obtain
\begin{equation}\label{ras62}
\sum_{n=0}^{\infty}\overline R_{6}^\ast (n)q^n\equiv F^2(q^3,q^6)+q^2\psi^2(q^9)+2q
F(q^3,q^6)\psi(q)\pmod 3.\end{equation}
Extracting the terms involving $q^{3n+2}$ from \eqref{ras62}, dividing by $q^2$ and replacing $q^3$ by $q$, we obtain 
\begin{equation}\label{ras62a}
\sum_{n=0}^{\infty}\overline R_{6}^\ast (3n+2)q^n\equiv\psi^2(q^3)\pmod 3.\end{equation}
Extracting the terms involving $q^{3n}$ from \eqref{ras62a} and replacing $q^3$ by $q$, we obtain
\begin{equation}\label{ras62b}
\sum_{n=0}^{\infty}\overline R_{6}^\ast (9n+2)q^n\equiv\psi^2(q)\pmod 3.\end{equation}
Employing \eqref{ras61a} in \eqref{ras62b}, we obtain
\begin{equation}\label{r1}
\sum_{n=0}^{\infty}\overline R_{6}^\ast (9n+2)q^n\equiv\sum_{n=0}^{\infty}\overline R_{6}^\ast (n)q^n \pmod 3.
\end{equation} Iterating \eqref{r1} with $n$ replaced by $9n+2$ and simplifying, we obtain
\begin{equation}\label{r1a}
\sum_{n=0}^{\infty}\overline R_{6}^\ast \left(9^\alpha n+\frac{9^\alpha-1}{2}\right)q^n\equiv\sum_{n=0}^{\infty}\overline R_{6}^\ast (n)q^n \pmod 3,
\end{equation} where $\alpha\ge0$ is any integer. Equating the coefficients of $q^n$ on both sides of \eqref{r1a}, we arrive at (i).\\
Combining \eqref{ras62} and \eqref{r1a}, we obtain
\begin{equation}\label{r1aw}
\sum_{n=0}^{\infty}\overline R_{6}^\ast \left(9^\alpha n+\frac{9^\alpha-1}{4}\right)q^n\equiv F^2(q^3,q^6)+q^2\psi^2(q^9)+2q
F(q^3,q^6)\psi(q)\pmod 3,\end{equation}  $\alpha\ge0$ is any integer.

Extracting the terms involving $q^{3n+2}$ from \eqref{r1aw}, dividing by $q^2$ and replacing $q^3$ by $q$, we obtain
\begin{equation}\label{r1aw1}
\sum_{n=0}^{\infty}\overline R_{6}^\ast \left(9^\alpha (3n+2)+\frac{9^\alpha-1}{4}\right)q^n\equiv \psi^2(q^3)\pmod 3,\end{equation} 
Extracting the terms involving $q^{3n+1}$ and $q^{3n+2}$ from \eqref{r1aw1}, we complete the proof of (ii) and (iii), respectively.
\end{proof}

\begin{theorem} Let $p\ge 3$  be a prime with $\left(\dfrac{-1}{p}\right)=-1$ and $1 \le r \le p-1$. Then for any integers $\alpha\geq0$ and $n\geq0$, we have
	$$\hspace{-5cm}(i)\quad \sum_{n=0}^{\infty}\overline R_{6}^\ast\left(2\cdot p^{2\alpha}n+\frac{5\cdot p^{2\alpha}-1}{4}\right) \equiv 2\psi(q)\psi(q^4) \pmod3,$$
	$$\hspace{-5.5cm}(ii)\quad \overline R_6^\ast\left(2\cdot p^{2\alpha+1}(pn+r)+\frac{5\cdot p^{2\alpha+2}-1}{4}\right) \equiv 0\pmod3.$$
\end{theorem}
\begin{proof} Employing  \eqref{b4} in \eqref{ras61}, we obtain
	\begin{equation}\label{rast6}
		\sum_{n=0}^{\infty}\overline R_{6}^\ast (n)q^n\equiv\frac{f_8^5}{f_2f_{16}^2}+2q\frac{f_4^2f_{16}^2}{f_2f_8} \pmod 3.
	\end{equation}
	Extracting the terms involving $q^{2n+1}$, dividing by $q$ and replacing $q^2$ by $q$ from \eqref{rast6}, we obtain
	\begin{equation}\label{rast6a}
		\sum_{n=0}^{\infty}\overline R_{6}^\ast (2n+1)q^n\equiv 2 \frac{f_2^2}{f_1}\cdot \frac{f_8^2}{f_4} \pmod{3}.
	\end{equation}
	Employing \eqref{tp1} in \eqref{rast6a}, we obtain
	\begin{equation}\label{rast6b}
		\sum_{n=0}^{\infty}\overline R_{6}^\ast (2n+1)q^n\equiv 2 \psi(q)\psi(q^4)\pmod{3}.
	\end{equation}
	which is $\alpha=0$ case of (i). Assume that (i) is true for some integer $\alpha\geq0$. Using \eqref{w4} in (i), we obtain
	$$\hspace{-9cm}\sum_{n=0}^{\infty}\overline R_{6}^\ast\left(2\cdot p^{2\alpha}n+\frac{5\cdot p^{2\alpha}-1}{4}\right) $$
     $$\hspace{-1cm}\equiv2\bigg[\sum_{j=0}^{(p-3)/2}{\it q}^{(j^2+j)/2} F\left( {\it q}^{(p^2+(2j+1)p)/2},{\it q}^{(p^2-(2j+1)p)/2}\right)+{\it q}^{(p^2-1)/8}\psi(q^{p^2})\bigg]$$
     \begin{equation}\label{rast6c}
	\hspace{1cm}\times\bigg[\sum_{t=0}^{(p-3)/2}{\it q}^{4(t^2+t)/2} F\left( {\it q}^{4(p^2+(2t+1)p)/2},{\it q}^{4(p^2-(2t+1)p)/2}\right)+{\it q}^{4(p^2-1)/8}\psi(q^{4p^2})\bigg]\pmod{3}.
    \end{equation}
	Consider the congruence
	$$\left(\frac{j^2+j}{2}\right)+4\left(\frac{t^2+t}{2}\right)\equiv \frac{5(p^2-1)}{8}\pmod{p} ,$$
	which is ĕquivalent to
	$$(2j+1)^2+(4t+2)^2\equiv0\pmod{p}.$$
	Since $\left(\dfrac{-1}{p}\right)=-1$, the above congruence has only solution $j=t=(\pm p-1)/2.$\\
Therefore, extracting the terms involving $q^{pn+5(p^2-1)/8}$ from both sides of \eqref{rast6c}, dividing throughout by $q^{5(p^2-1)/8}$ and then replacing $q^p$ by $q$, we obtain
		\begin{equation}\label{rast6d}
			\sum_{n=0}^{\infty}\overline R_{6}^\ast\left(2\cdot p^{2\alpha+1}n+\frac{5\cdot p^{2\alpha+2}-1}{4}\right) \equiv 2\psi(q^p)\psi(q^{4p}) \pmod3.
		\end{equation}
			Extracting the terms involving $q^{pn}$ from \eqref{rast6d} and replacing $q^p$ by $q$, we obtain
			\begin{equation}\label{rast6e}
				\sum_{n=0}^{\infty}\overline R_{6}^\ast\left(2\cdot p^{2(\alpha+1)}n+\frac{5\cdot p^{2(\alpha+1)}-1}{4}\right) \equiv 2\psi(q)\psi(q^4) \pmod3
			\end{equation}
			which is $\alpha+1$ case of (i). Thus, by  induction, we complete the proof of (i).
			Extracting the terms involving $q^{pn+r}$, for $1\le r\le p-1$, from both sides \eqref{rast6d}, we complete the proof (ii).
			\end{proof}

\begin{theorem}We have
$$\hspace{-8.6cm}(i)\quad\overline R_{8}^\ast (4n+\xi)\equiv 0\pmod 4; \quad \xi=2,3,$$
$$\hspace{-7cm}(ii)\quad\overline R_{8}^\ast (16n+4\xi+1)\equiv 0\pmod 4,\quad \xi=1,2,3,$$
\end{theorem}
\begin{proof} Setting $\ell=8$ in \eqref{rast}, we obtain
	\begin{equation}\label{rast8}
		\sum_{n=0}^{\infty}\overline R_{8}^\ast (n)q^n=\frac{f_2f_8}{f_1^2}.
	\end{equation}
	Employing \eqref{b4} in \eqref{rast8}, we obtain
	\begin{equation}\label{rast8a}
		\sum_{n=0}^{\infty}\overline R_{8}^\ast (n)q^n=\frac{f_8^6}{f_2^4f_{16}^2}+2q\frac{f_4^2f_{16}^2}{f_2^4}.
	\end{equation}
 Employing \eqref{yp2} with $p=2$ in \eqref{rast8a}, we obtain 
	\begin{equation}\label{rast8b}
			\sum_{n=0}^{\infty}\overline R_{8}^\ast (n)q^n\equiv \frac{f_8^2}{f_4^2}+2qf_{16}^2\pmod{4}.
	\end{equation}
	Extracting the terms involving $q^{4n+\xi}$ for $\xi=2,3$ from \eqref{rast8b}, we arrive at (i).\\ Again, extracting the terms $q^{4n+1}$ from \eqref{rast8b}, dividing by $q$ and replacing $q^4$ by $q$, we obtain
	\begin{equation}\label{rast8c}
			\sum_{n=0}^{\infty}\overline R_{8}^\ast (4n+1)q^n\equiv 2f_4^2\pmod{4}.
	\end{equation}
	Extracting the terms involving $q^{4n+\xi}$ for $\xi=1,2,3$ from \eqref{rast8c}, we arrive at (ii).\\ 
\end{proof}
\begin{theorem}We have
	$$\hspace{-10.4cm}(i)\quad\overline R_{8}^\ast (4n+3)\equiv 0\pmod 8,$$
	$$\hspace{-7cm}(ii)\quad\overline R_{8}^\ast (8n+2\xi+1)\equiv 0\pmod 8,\quad \xi=1,2,3.$$
\end{theorem}
\begin{proof}Extracting the terms involving $q^{2n+1}$ from \eqref{rast8a}, dividing by $q$ and replacing $q^2$ by $q$, we obtain
	\begin{equation}\label{ra8}
		\sum_{n=0}^{\infty}\overline R_{8}^\ast (2n+1)q^n=2\frac{f_2^2f_8^2}{f_1^4}.
	\end{equation}
	Employing \eqref{b3} in \eqref{ra8}, we obtain
	\begin{equation}\label{ra81}
		\sum_{n=0}^{\infty}\overline R_{8}^\ast (2n+1)q^n=2\frac{f_4^{14}}{f_2^{12}f_8^2}+8q\frac{f_4^2f_8^6}{f_2^8}.
	\end{equation}
	Employing \eqref{yp2} with $p=2$ in \eqref{ra81}, we obtain
	\begin{equation}\label{r8}
		\sum_{n=0}^{\infty}\overline R_{8}^\ast (2n+1)q^n\equiv 2f_8^2\pmod{8}
	\end{equation}
	Extracting the terms involving $q^{2n+1}$ on both sides of \eqref{r8}, we arrive at (i).\\
  Similarly, (ii) follows from  \eqref{r8} by extracting the terms involving $q^{4n+\xi}$ for $\xi=1,2,3$.
\end{proof}

\begin{theorem} Let $p\ge 5$  be a prime with $\left(\dfrac{-3}{p}\right)=-1$ and $1 \le r \le p-1$. Then for any integers $\alpha\geq0$ and $n\geq0$, we have
    $$\hspace{-5cm} (i) \quad \sum_{n=0}^{\infty}\overline R_{8}^\ast\left(8p^{2\alpha}n+\frac{4\cdot p^{2\alpha}-1}{3}\right)q^n\equiv2f_1\psi(q)\pmod{8},$$
	 
	$$\hspace{-5.3cm} (ii) \quad \overline R_{8}^\ast\left(8p^{2\alpha+1}(pn+r)+\frac{4\cdot p^{2\alpha+2}-1}{3}\right)\equiv 0\pmod{8},$$
	$$\hspace{-4cm}(iii) \quad \frac{1}{2}~\overline R_8^\ast (~16n+1)\equiv p(n)+\sum_{\nu=1}^{\infty} p\left(n-\frac{\nu(\nu+1)}{2}\right)\pmod 4.$$
\end{theorem}
\begin{proof}Extracting the terms involving $q^{4n}$ from \eqref{r8} and simplifying, we obtain
	\begin{equation}\label{rast8d}
		\sum_{n=0}^{\infty}\overline R_{8}^\ast (8n+1)q^n\equiv2\frac{f_1f_2^2}{f_1}\pmod{8}.
	\end{equation}
	Employing \eqref{tp1} in \eqref{rast8d},we obtain
	\begin{equation}\label{rast8e}
		\sum_{n=0}^{\infty}\overline R_{8}^\ast (8n+1)q^n\equiv2f_1\psi(q)\pmod{8}.
	\end{equation}
	which is the $\alpha=0$ case of (i). Assume that (i) is true for some integer $\alpha\geq0$. Using \eqref{w4} in \eqref{u10}, we find that
	$$\hspace{-9cm}\sum_{n=0}^{\infty}\overline R_{8}^\ast\left(8p^{2\alpha}n+\frac{4\cdot p^{2\alpha}-1}{3}\right)q^n$$
	$$\equiv2\bigg[\sum_{\substack{ k={-(p-1)/2} \\ k \ne {(\pm p-1)/6}}}^{(p-1)/2}(-1)^{k}{q}^{\frac{3k^2+k}{2}} F\left( {-q}^{(3p^2+(6k+1)p)/2},{-q}^{(3p^2-(6k+1)p)/2}\right)+(-1)^{{(\pm p-1)/6}}{q}^{(p^2-1)/24}f_{p^2}\bigg]$$
	\begin{equation}\label{rast8f}
		\hspace{.5cm}\times\bigg[\sum_{j=0}^{(p-3)/2}{\it q}^{(j^2+j)/2} F\left( {\it q}^{(p^2+(2j+1)p)/2},{\it q}^{(p^2-(2j+1)p)/2}\right)+{\it q}^{(p^2-1)/8}\psi(q^{p^2})\bigg]\pmod{8}.
	\end{equation}
	Consider the congruence
	$$ \frac{(3k^2+k)}{2}+\frac{(j^2+j)}{2}\equiv\frac{4(p^2-1)}{24}\pmod{p},$$
	which is equivalent to 
	$$ (6k+1)^2+3(2j+1)^2\equiv0\pmod{p}.$$
	Since $\left(\dfrac{-3}{p}\right)=-1$, the above congruence has only solution $k=(\pm p-1)/6$ and $j=(\pm p-1)/2$.
	Therefore, extracting the terms involving $q^{pn+(p^2-1)/6}$ from both sides of \eqref{rast8f}, dividing throughout by $q^{(p^2-1)/6}$ and then replacing $q^p$ by $q$, we obtain
	\begin{equation}\label{rast8g}
		\sum_{n=0}^{\infty}\overline R_{8}^\ast\left(8p^{2\alpha+1}n+\frac{4p^{2\alpha+2}-1}{3}\right)q^n\equiv2f_p\psi(q^p)\pmod{8}.
	\end{equation}
	Extracting the terms involving $q^{pn}$ from \eqref{rast8g} and replacing $q^p$ by $q$, we obtain
	\begin{equation}\label{rast8h}
		\sum_{n=0}^{\infty}\overline R_{8}^\ast\left(8p^{2(\alpha+1)}n+\frac{4p^{2\alpha+2}-1}{3}\right)q^n\equiv2f_1\psi(q)\pmod{8},
	\end{equation}
	which is $\alpha+1$ case of (i). Thus, by induction the proof of (i) complete.
	
	Extracting the  terms involving $q^{pn+r}$, for $1\le r\le p-1$, from both sides \eqref{rast8g}, we  arrive at (ii).

	Extracting the terms involving $q^{16n+1}$ from \eqref{rast8b}, dividing by $q$ and replacing $q^4$ by $q$, we obtain 
	 		\begin{equation}\label{ras8c}
		\sum_{n=0}^{\infty}\overline R_{8}^\ast (16n+1)q^n\equiv 2f_1^2\equiv2\frac{f_2^2}{f_1^2}\pmod{4}, 
	\end{equation}where we used \eqref{yp2} with $p=2$.

	Employing \eqref{pngeni} and \eqref{tp1} in \eqref{ras8c}, we obtain
	\begin{equation}\label{ras8d}
		\frac{1}{2}\sum_{n=0}^{\infty}\overline R_{8}^\ast (16n+1)q^n\equiv\sum_{n=0}^{\infty}p(n)q^n+\sum_{\nu=0}^{\infty}p\left(n-\nu(\nu+1)/2\right)q^n \pmod{4}.
	\end{equation}
	Equating the coefficient of $q^n$ on both sides of \eqref{ras8d}, we arrived at (iii).
\end{proof}

\begin{theorem}\ We have
	$$\hspace{-8.4cm}(i)\quad\overline R_\ell^\ast (n)+\sum_{\nu=1}^{\infty}\overline R_\ell^\ast (n-\ell\nu)p(\nu)=\overline p(n),$$
	$$\hspace{-11.8cm}(ii)\quad\overline R_2^\ast (n)=D_2(n).$$where $D_2(n)$ is the number of partitions of $n$ into distinct parts such that each part has two copies.\end{theorem}
\begin{proof}
From \eqref{rast}, we have
	\begin{equation}\label{l}
		\sum_{n=0}^{\infty}\overline R_\ell ^\ast (n)q^n=\frac{f_2f_\ell}{f_1^2}
	\end{equation}
		Using \eqref{opn} in \eqref{l}, we obtain
	\begin{equation}\label{l1}
			\sum_{n=0}^{\infty}\overline R_\ell ^\ast (n)q^n\cdot\frac{1}{f_\ell}=\sum_{n=0}^{\infty}\bar{p}(n)q^n.
	\end{equation}Employing \eqref{pngeni} in \eqref{l1}
	\begin{equation}\label{l2}
		\sum_{n=0}^{\infty}\overline R_\ell ^\ast (n)q^n\left(1+\sum_{\nu=1}^{\infty}p(\nu)q^{\ell\nu}\right)=\sum_{n=0}^{\infty}\bar{p}(n)q^n
	\end{equation}
	Extracting the coefficient of $q^n$ from both sides of \eqref{l2}, we arrive at (i).

Setting $\ell=2$ in \eqref{l} and simplifying,  we obtain
\begin{equation}\label{lr}
		\sum_{n=0}^{\infty}\overline R_\ell ^\ast (n)=\frac{f_2^2}{f_1^2}=(-q;q)_\infty^2.
	\end{equation} Now (ii) follows immediately from \eqref{lr}.
\end{proof}

\section*{\bf Acknowledgement}  The second author thanks University Grants Commission (UGC) of India for supporting her research work through CSIR-UGC Junior Research Fellowship (JRF) vide  NTA-Ref. No.: 221610074648, Dated 29/11/2022.

\section*{\bf Declarations}

\noindent{\bf	Author Contributions.} Both authors contributed equally in the article. 

\noindent{\bf Conflict of Interest.} The authors declare that there is no conflict of interest regarding the publication of
this article.

\noindent{\bf Human and animal rights.} The authors declare that there is no research involving human participants or
animals in the contained of this paper.	

\noindent{\bf Data availability statements.} Data sharing not applicable to this article as no datasets were generated or analysed during the current study.


\begin{thebibliography}{99}
   \bibitem{alanzi} Alanazi A. M.,  Alenazi B. M., W. J. Keith and  Munagi A. O.: Refining overpartitions by properties of non overlined parts, \textit{Contrib. Discrete Math.} \textbf{17}(2) (2022), 96-111.
   
   	\bibitem{nak} Ballantine, C. and Merca, M.:\newblock { Almost $3$-regular overpartitions}, \textit{Ramanujan J.} \textbf{58} (2022), 957-971.
   
   \bibitem{zlndb}  Baruah, N. D. and Ahmed, Z.: New congruences for $\ell$-regular partitions for $\ell\in \{5, 6, 7, 49\}$, {\it The Ramanujan J.} {40} (2016), 649-668.
   
  	\bibitem{klndb}  Baruah, N. D. and Das, K.: Parity results for $7$-regular and $23$-regular partitions, {\it Int. J.  Number Theory} {\bf 11} (2015), 2221-2238.
	
	 
	\bibitem{bcb3} Berndt,  B. C.: {\it Ramanujan's Notebooks, Part	III}. Springer-Verlag, New York, 1991. 
	
	\bibitem{bc1} Berndt, B. C.  and Rankin, R. A.: \textit{Ramanujan: Letters and Commentary}, Amer. Math. Soc. (1995).
	
	\bibitem{ps} Buragohain, P. and Saikia, N.: New congruences for overpartitions with $\ell$-regular overlined parts.\textit{J. Anal.} \textbf{31} (2023), 1819-1837.
	
	\bibitem{carlson} Carlson, R. and Webb, J. J.: Infinite families of congruences for $k$-regular partitions, {\it Ramanujan J.} {\bf 33} (2014), 329-337.
	
	\bibitem{cor} Corteel, S. and Lovejoy, J.: {Overpartitions}, {\it  T. Am. Math. Soc. Comp.} \textbf{356}, (2004), 1623-1635.
	
	\bibitem{cui} Cui, S. P.  and Gu, N. S. S.: Arithmetic properties of $l$-regular partitions, {\it Adv. Appl. Math.} {\bf 51} (2013) 507-523.
	
	\bibitem{penniston}  Dandurand, B. and Penniston, D.: ${\ell}$-divisibility of ${\ell}$-regular partition functions, {\it Ramanujan J.} {\bf 19} (2009), 63-70.
	
	\bibitem{ovyt} Hirschhorn, M. D.: {\it The Power of $q$,} A Personal Journey, Developments in Mathematics, {\bf 49.} Springer, Cham (2017).
	
	\bibitem{hs2} Hirschhorn, M. D. and Sellers, J. A.: An infinite family of overpartition congruences modulo 12, {\it Integers} {\bf 5} (2005), A20 1-4. 
	
	\bibitem{hs1} Hirschhorn, M. D. and Sellers, J. A.: Arithmetic relations for overpartitions, {\it J. Combin. Math. Combin. Comput.} {\bf  53} (2005), 65-73.
	
	\bibitem{hirs}  Hirschhorn, M. D. and Sellers, J. A.: Elementary proofs of parity results for $5$-regular partitions, {\it Bull. Aust. Math. Soc.} {\bf 81} (2010), 58-63.
	
	\bibitem{kim1} Kim, B.: A short note on the overpartition function, {\it Discrete Math.} {\bf 309} (2009), 2528-2532.
	
	\bibitem{lin2}  Lin, B. L. S.: An infinite family of congruences modulo $3$ for $13$-regular bipartitions, {\it Ramanujan J.} {\bf 39}(1) (2016), 169-178.
	
	\bibitem{kmah}  Mahlburg, K.: The overpartition function modulo small powers of $2$, {\it Discrete Math.}	{\bf286} (2004), 263-267.
	
	\bibitem{penniston1} Penniston, D.: Arithmetic of ${\ell}$-regular partition functions, {\it Int. J. Number Theory} {\bf 4} (2008), 295-302.
	
	\bibitem{sr2} Ramanujan, S.: Congruence properties of partitions, {\it Math. Z.} {\bf9} (1921), 147-153.
	
	\bibitem{SR} Ramanujan, S. and Hardy, G. H.:\textit{Collected papers}, Chelsea, New York, (1962).
	
	\bibitem{ij} Saikia, N. and Boruah, C.: Congruences for $\ell$-regular overpartition for $\ell\in\{5, 6, 8\}$.\textit{ Indian J. Pure Appl. Math.} \textbf{48}(2) (2017), 295-308.
	
	\bibitem{sq}  Sellers, J. A.: Extending congruences for overpartitions with $\ell$-regular nonoverlined parts, {\it Bull. Aust. Math. Soc.} (2024), doi:10.1017/S0004972724001023.
	
	\bibitem{sc} Shivanna, G. D.  and Chandrappa, S.: Congruences for overpartitions with $l$-regular overlined parts, {\it J. Anal.} \textbf{31}  (2023), 1819-1837
	
	\bibitem{webb}  Webb, J. J.: Arithmetic of the $13$-regular partition function modulo 3, {\it Ramanujan J.} {\bf 25} (2011), 49-56.
	
	 \bibitem{xia1over} Xia, E. X. W.: Congruences modulo $9$ and $27$ for overpartitions, {\it Ramanujan J.} {\bf 42}(2) (2017), 301-323.
	
	\bibitem{xia1} Xia, E. X. W. and Yao, O. X. M.: A proof of Keith's conjecture for $9$-regular partitions modulo $3$, {\it Int. J. Number Theory} {\bf 10} (2014), 669-674.
  \end{thebibliography}
 \end{document}